\documentclass[12pt]{article}

\makeatletter
\renewcommand{\@oddfoot}{\hfill \thepage}
\makeatother

\usepackage[english]{babel}
\usepackage{mathrsfs}
\usepackage{amssymb}
\usepackage{amsmath,amsthm,amsfonts}
\usepackage{cite}
\usepackage[dvipdf]{graphicx}
\usepackage{fancyhdr}
\usepackage[titletoc]{appendix}
\usepackage{epstopdf} 

\begin{document}

\begin{center}
\Large\bf{A diffusion model of surface soil pollution based on planar finite-velocity 
            stochastic motion\\ with random lifetime}
\end{center}

\begin{center}
Alexander D. KOLESNIK\\
Institute of Mathematics and Computer Science\\ 
Moldova State University\\
Academy Street 5, Kishinev 2028, Moldova
\end{center}

{\bf Abstract} 
We present a diffusion model of surface soil pollution from a stationary source based on the symmetric  stochastic motion at finite speed in the plane $\Bbb R^2$, also called the planar Markov random flight, whose lifetime is a random variable with given distribution. We consider a heavy-particle model, in which 
the lifetime is supposed to be an exponentially-distributed random variable, and obtain an explicit formula  for the stationary probability density of the pollution process expressed in terms of McDonald functions with variable indices. We also study a light-particle model, in which the lifetime is a gamma-distributed random variable. In this case, the stationary probability density of the pollution process is given in the form of a definite integral calculated numerically, as well as in the form of a functional series composed of the hypergeometric functions with variable coefficients. These stationary densities are plotted in a figure and numerically calculated tables that demonstrate the behaviour of the pollution process on long time intervals. Some remarks on the pollution model based on asymmetric finite-velocity planar stochastic motion are also given. 

\bigskip  


\bigskip 


\bigskip 

\section{Introduction} 
\setcounter{equation}{0} 

Mathematical models of environmental pollution have become an important field of modern researches in applied mathematics and ecology. Besides various models of air and water pollution, the processes of soil pollution are of a special importance. For such processes, two main classes of models should be 
distinguished. The first type of such models is related to the problem of describing the process of 
surface soil pollution from a stationary source. The second one is related to the problem of describing 
the process of seepage of a pollutant from the surface to the deep strata of soil. It is obvious that these two classes of problems are connected with each other, since the characteristics of surface soil pollution  can serve as initial and/or boundary conditions for describing the process of pollutant percolation into its deeper strata. The reader interested in some modern mathematical models of soil pollution can address to the article \cite{stag} and bibliography therein. 

The main mathematical tools in constructing such models are differential equations, both ordinary and partial. This approach implies the attempts to take into account several basic parameters depending on time and to describe their behaviour and connections with each other by means of a system of nonlinear differential and integro-differential equations. It usually consists of sometimes fairly big number of nonlinear equations with non-trivial initial and/or boundary conditions and solving of such a system, even numerically, is an extremely complicated mathematical problem. 

Another approach implies the use of stochastic methods for describing various pollution processes. The main tools in studying such models are the stochastic partial differential equations (see, for instance, 
\cite{kalli} and bibliography therein, where several stochastic models arising from environmental problems are considered).  

In this article we present a diffusion model of surface soil pollution from a stationary source based on 
the stochastic motion at finite speed in the plane, also called {\it the planar Markov random flight}, 
whose lifetime is a positive random variable with given distribution. The main feature of our approach is the refusal of attempts to take into account the influence of several factors related to each other by various functional relationships and leading to complex systems of nonlinear equations, both deterministic and stochastic. Instead, we propose to consider not the influence of such factors separately, but only the result of their interaction, that is, the level of pollution. With this approach, it is natural to assume that pollution is carried out by some particles, which after leaving the source perform a random walk with a finite speed and the time of their active evolution (lifetime) is a positive random variable with a known distribution. With this interpretation, the pollution level is nothing more than the density of the location of polluting particles on the plane. Moreover, it is obvious that it makes sense to consider the pollution density only over long time intervals, that is, the stationary density and this is the main subject of our interest.  

The article is organized as follows. In Section 2, for the reader's, convenience, we give a succinct description of the planar symmetric Markov random flight, the structure of its distribution, as well as the explicit formula for its transition density, which our model is based on. In Section 3 we describe a  general model of surface soil pollution from a stationary source represented by the planar symmetric Markov random flight with random lifetime and the stationary density of the pollution process is given in the form of a definite integral. A heavy-particle model related to the case when the lifetime is an exponentially distributed random variable is considered in Section 4. The stationary density of soil pollution is obtained in the form of a series with respect to McDonald functions with variable indices. A numerical example is also given. A similar analysis is presented in Section 5 for the light-particle model related to the case when the lifetime is a gamma-distributed random variable. The stationary density of soil pollution is given  in the form of a definite integral, as well as in the form of a series with respect to hypergeometric functions. A  numerical example is also considered. In Section 6, we present some remarks related to the  asymmetric model based on the planar stochastic motion, when the directions are taken on according to the von Mises distribution. A general formula for the stationary density of soil pollution in the form of a series of definite integrals whose integrands are the double convolutions with respect to spatial and time variables, is presented.

\section{Planar symmetric finite-velocity stochastic motion} 
\setcounter{equation}{0} 

In this section we describe the planar symmetric finite-velocity stochastic motion, which our  
diffusion model is based on, and give some known results related to the structure and closed-form 
expression of its distribution. 

Consider the stochastic motion in the Euclidean plane $\Bbb R^2$   
represented by a particle that, at the initial time instant $t=0$, starts from the 
origin $\bold 0 = (0, 0)\in\Bbb R^2$ and moves with some finite speed $c$ 
(note that $c$ is treated as the constant norm of the velocity). The initial direction is a random 
two-dimensional unit vector uniformly distributed on the unit circumference 
$S_1 = \left\{ \bold x=(x_1, x_2)\in \Bbb R^2: \;
\Vert\bold x\Vert^2 = x_1^2+x_2^2=1 \right\}$ .

The particle changes its direction at random time instants that form a homogeneous Poisson 
flow of rate $\lambda>0$. In each of such Poissonian moments the particle instantly takes on 
a new random direction uniformly distributed on $S_1$ independently of its previous direction. 
Each sample path of this motion represents a planar broken line of total length $ct$ composed of 
the segments of $\lambda$-exponentially distributed random lengths and uniformly oriented 
in $\Bbb R^2$. The trajectories of the process are continuous and differentiable almost everywhere.

Let $\bold X(t)=(X_1(t), X_2(t))$ be the particle's position at 
arbitrary time instant $t>0$. The stochastic process $\bold X(t)$ is referred to as the 
{\it planar symmetric Markov random flight}. Let $d\bold x$ denote an infinitesimal 
element in the plane $\Bbb R^2$ with Lebesgue measure $\nu(d\bold x) = dx_1 dx_2$.

At arbitrary time instant $t>0$ the particle, with probability 1, is located in the 
circle of radius $ct$: 
\begin{equation}\label{3.1.2}
\bold B_{ct} = \left\{ \bold x=(x_1, x_2)\in \Bbb R^2 : \;
\Vert\bold x\Vert^2 = x_1^2+x_2^2\le c^2t^2 \right\} .
\end{equation}

The distribution 
\begin{equation}\label{3.1.3}
\text{Pr} \left\{ \bold X(t)\in d\bold x \right\} = \text{Pr}
\left\{ X_1(t)\in dx_1, X_2(t)\in dx_2 \right\} , \qquad \bold x\in\bold B_{ct}, \quad t>0, 
\end{equation}
consists of two components. 

The singular component of distribution (\ref{3.1.3}) 
is related to the case when no Poisson events occur in the time interval $(0,t)$ 
and, therefore, the particle does not change its initial direction. The singular component 
is concentrated on the circumference of radius $ct$: 
\begin{equation}\label{3.1.4}
S_{ct} =\partial\bold B_{ct} = \left\{ \bold x=(x_1, x_2)\in \Bbb R^2: \; 
\Vert\bold x\Vert^2 = x_1^2+x_2^2=c^2t^2 \right\} . 
\end{equation}

If at least one Poisson event occurs in the time interval $(0,t)$ and, therefore, 
the particle at least once changes its direction, then, at time $t$, it is located 
strictly inside the circle $\bold B_{ct}$ and the part of distribution (\ref{3.1.3}) 
corresponding to this case is concentrated in the interior of the circle $\bold B_{ct}$: 
\begin{equation}\label{3.1.7}
\text{int} \; \bold B_{ct} = \left\{ \bold x=(x_1,x_2)\in \Bbb R^2: \; 
\Vert\bold x\Vert^2 = x_1^2+x_2^2<c^2t^2 \right\} 
\end{equation}
and forms its absolutely continuous component. 

Therefore, there exists the density of the absolutely continuous component of the distribution (\ref{3.1.3})
\begin{equation}\label{3.1.8}
p^{(ac)}(\bold x,t) = f(\bold x,t) \Theta(ct-\Vert\bold x\Vert) , \qquad 
\bold x\in \text{int} \; \bold B_{ct} , \quad t>0,
\end{equation}
where $f(\bold x,t)$ is some positive function absolutely continuous in $\text{int} \; \bold B_{ct}$ 
and $\Theta(x)$ is the Heaviside unit-step function. The existence of 
density (\ref{3.1.8}) follows from the fact that, since the sample paths of the process $\bold X(t)$ 
are continuous and differentiable almost everywhere, the distribution (\ref{3.1.3}) must contain an absolutely continuous component and this justifies the existence of density (\ref{3.1.8}). 

Hence, the density of distribution (\ref{3.1.3}) has the structure
$$p(\bold x,t) = p^{(s)}(\bold x,t) + p^{(ac)}(\bold x,t), 
\qquad \bold x\in\bold B_{ct} , \quad t>0,$$
where $p^{(s)}(\bold x,t)$ and $p^{(ac)}(\bold x,t)$ are the densities (in the sense 
of generalized functions) of the singular and absolutely continuous components 
of distribution (\ref{3.1.3}), respectively. 

The most important result is represented by the density of distribution (\ref{3.1.3}) obtained,  
by different methods, in a series of works 
(see \cite[formula (5.2.5)]{kol1},\cite{kol4},\cite[formula (21)]{kol5},\cite{mas},\cite{sta2}) 
and given by the relation : 
\begin{equation}\label{3.1.9}
p(\bold x,t) = \frac{e^{-\lambda t}}{2\pi ct} \; \delta(c^2t^2 - \Vert\bold x\Vert^2) + 
\frac{\lambda}{2\pi c} \; \frac{\exp \left( -\lambda t +
\frac{\lambda}{c} \sqrt{c^2t^2 - \Vert\bold x\Vert^2} \right)}{\sqrt{c^2t^2-\Vert\bold x\Vert^2}} \; 
\Theta(ct-\Vert\bold x\Vert) , 
\end{equation}
$$\bold x = (x_1, x_2)\in\Bbb R^2, \qquad \Vert\bold x\Vert^2 = x_1^2+x_2^2, \qquad t>0 .$$

The first term in (\ref{3.1.9}) is the singular part of density concentrated on the circumference 
$S_{ct}$, while its second term 
\begin{equation}\label{AbsCont}
p^{(ac)}(\bold x,t) = \frac{\lambda}{2\pi c} \; \frac{\exp \left( -\lambda t +
\frac{\lambda}{c} \sqrt{c^2t^2 - \Vert\bold x\Vert^2} \right)}{\sqrt{c^2t^2-\Vert\bold x\Vert^2}} \; 
\Theta(ct-\Vert\bold x\Vert) , 
\end{equation}
represents the absolutely continuous part of density (\ref{3.1.9}) concentrated in 
$\text{int} \; \bold B_{ct}$ and plays important role in our diffusion model.

\section{General diffusion model of surface soil pollution} 

In this section we present a general diffusion model of surface soil 
pollution from a stationary source based on the planar symmetric 
Markov random flight described in previous section. 

The soil surface pollution can be imagined as follows. A stationary source (a pipe of an industrial enterprise, for example) emits polluting particles that settle on the soil. The particles have random masses and, at the moment of emitting of each particle, the strength and direction of the wind are also random. When flying, the particle makes chaotic three-dimensional movements that are random both on their directions and the length of free runs. All these random factors determine a place where the particle falls on the soil. The time of active evolution of the particle, that is, the time from the moment of its emitting to the moment it falls onto the soil, we call the {\it particle's lifetime} and this is a positive random variable. The degree of pollution of a certain part of the territory around the source at some given time moment is directly proportional to the number of particles that have fallen to this site by this time. 

Since we are interested not in the three-dimensional motion itself after the particle leaves the source, but in the distribution of the place of its fall to the ground, it is natural to focus on studying the behaviour of the projection of three-dimensional random motion onto a plane. Our {\it basic assumption} is that the projection of a three-dimensional random motion of a particle, on large time intervals, is asymptotically a two-dimensional symmetric Markov random flight, whose lifetime is a positive random variable with given distribution. Such assumption is somewhat justified because, as it was shown in \cite{kol2}, the similar connection exists between the marginals of the planar symmetric Markov random flight and the one-dimensional Goldstein-Kac telegraph process whose densities are asymptotically equivalent, as time $t\to\infty$.  
Such interpretation allows us to apply the known results related to the planar Markov random flight and  presented in previous section, for obtaining the stationary probability density of the pollution process. 

Let $q(t)$ be the density of the distribution of the random particle's lifetime, assumed to exist. 
Then, according to (\ref{AbsCont}), the stationary density of the distance $\Vert\bold x\Vert$ between 
the fallen particle and emitting source is given by the general formula: 
\begin{equation}\label{PollutionDensityGeneral}
\aligned 
p(\bold x) & = \frac{\lambda}{2\pi c} \int_0^{\infty} \frac{\exp \left( -\lambda\tau +
\frac{\lambda}{c} \sqrt{c^2\tau^2 - \Vert\bold x\Vert^2} \right)}{\sqrt{c^2\tau^2-\Vert\bold x\Vert^2}} \; 
\Theta(c\tau-\Vert\bold x\Vert) \; q(\tau) \; d\tau \\ 
& = \frac{\lambda}{2\pi c} \int_{\Vert\bold x\Vert/c}^{\infty} e^{-\lambda\tau} \; 
\frac{\exp \left(\frac{\lambda}{c} \sqrt{c^2\tau^2 - \Vert\bold x\Vert^2} \right)}{\sqrt{c^2\tau^2 
- \Vert\bold x\Vert^2}} \; q(\tau) \; d\tau . 
\endaligned
\end{equation}
Note that, since the singular part of density (\ref{3.1.9}) tends to zero very quickly, as $t\to\infty$,
we may ignore it and take into account only the absolutely continuous part of density (\ref{3.1.9}). 

Formula (\ref{PollutionDensityGeneral}) represents the general relation for evaluating the density of pollution at the distance $\Vert\bold x\Vert$ from the source (origin $\bold 0 = (0,0)\in\Bbb R^2$) for large values of time $t$, that is, the stationary density of pollution at a point $\bold x = (x_1, x_2)\in\Bbb R^2$. The density $q(t)$ of the random lifetime of a particle can be selected in each particular case proceeded from the specifics of the model. 

In the next sections we will consider two particular cases of our general diffusion model when the emitted particles can have heavy and light weights.

\section{Heavy-particle diffusion model} 

Suppose that the source emits heavy-weighted polluting particles. In this case, it is natural to assume 
that most of the particles fall near the source and only a small portion of them may fall at a greater distance. Therefore, we may imagine the lifetime as a positive-valued random variable whose density has maximum at the origin $t=0$ and is monotonously decreasing, as time $t$ increases. In this case, the most appropriate is a $\mu$-exponentially distributed random variable with the density: 

\begin{equation}\label{LifetimeDensity1}
q(t) = \left\{ 
\aligned  \mu \; e^{-\mu t}, \quad & \text{if} \; t\ge 0, \\ 
          0, \qquad \; & \text{if} \; t<0 ,  
\endaligned\right.
\qquad \mu>0. 
\end{equation}

Then, according to general formula (\ref{PollutionDensityGeneral}), the stationary density of pollution 
at the distance $\Vert\bold x\Vert$ from the source is given by the formula: 
 
\begin{equation}\label{PollutionDensity1}
p_h(\bold x)  = \frac{\lambda\mu}{2\pi c} \int_{\Vert\bold x\Vert/c}^{\infty} e^{-(\lambda+\mu)\tau} \; 
\frac{\exp \left(\frac{\lambda}{c} \sqrt{c^2\tau^2 - \Vert\bold x\Vert^2} \right)}{\sqrt{c^2\tau^2 
- \Vert\bold x\Vert^2}} \; d\tau , 
\end{equation}
where, we remind the reader, $\lambda$ and $\mu$ are the parameters of the exponential distributions of the intensity of changes of direction and the particle's lifetime, respectively. 

In order to obtain the stationary density, we need to evaluate (\ref{PollutionDensity1}). This result is given by the following theorem. 

\bigskip 

{\bf Theorem 1.} {\it The stationary probability density $p_h(\bold x)$ of the symmetric 
planar Markov random flight with $\mu$-exponentially distributed random lifetime 
(\ref{LifetimeDensity1}), is given by the formula:} 

\begin{equation}\label{eq2.0}
p_h(\bold x) = \frac{\lambda\mu}{2\pi^{3/2} \; c^2} \; \sum_{k=0}^{\infty} \frac{\lambda^k}{k!} \;  
\Gamma\left( \frac{k+1}{2} \right) \; \left( \frac{2\Vert\bold x\Vert}{c(\lambda+\mu)} \right)^{k/2} \; 
K_{k/2}\left( (\lambda+\mu) \frac{\Vert\bold x\Vert}{c} \right) ,
\end{equation}
{\it where $K_{k/2}(z)$ is the McDonald function of order $k/2$.}

\begin{proof} 

By expanding the exponential function in the fraction of the integrand in (\ref{PollutionDensity1}), 
we obtain: 
$$\aligned
p_h(\bold x) & = \frac{\lambda \mu}{2\pi c} \; \int_{\Vert\bold x\Vert/c}^{\infty} \; 
e^{-(\lambda+\mu) \tau} \; \; \frac{\exp\left(\frac{\lambda}{c} \sqrt{c^2\tau^2 -
\Vert\bold x\Vert^2} \right) }{\sqrt{c^2\tau^2 - \Vert\bold x\Vert^2}} \; d\tau \\ 
& = \frac{\lambda\mu}{2\pi c} \int_{\Vert\bold x\Vert/c}^{\infty} e^{-(\lambda+\mu) \tau} \; 
\sum_{k=0}^{\infty} \frac{1}{k!} \; \left( \frac{\lambda}{c} \right)^k \; \left( \sqrt{c^2\tau^2 - 
\Vert\bold x\Vert^2} \right)^{k-1} \; d\tau \\ 
& = \frac{\lambda\mu}{2\pi c^2} \sum_{k=0}^{\infty} \frac{\lambda^k}{k!} \; 
\int_{\Vert\bold x\Vert/c}^{\infty} e^{-(\lambda+\mu) \tau} \; \left( \tau^2 - 
\frac{\Vert\bold x\Vert^2}{c^2} \right)^{(k-1)/2} \; d\tau \\ 
& \qquad (\text{see \cite[item 2.3.5, formula 4]{pbm2}}) \\ 
& = \frac{\lambda\mu}{2\pi^{3/2} \; c^2} \; \sum_{k=0}^{\infty} \frac{\lambda^k}{k!} \;  
\Gamma\left( \frac{k+1}{2} \right) \; \left( \frac{2\Vert\bold x\Vert}{c(\lambda+\mu)} \right)^{k/2} \; 
K_{k/2}\left( (\lambda+\mu) \frac{\Vert\bold x\Vert}{c} \right) ,
\endaligned$$
where $K_{k/2}(z)$ is the McDonald function of order $k/2$, also called the modified Bessel function of 
second kind. 

The theorem is proved.  
\end{proof}

The shape of density (\ref{eq2.0}) is plotted in Fig. 1 below.
\begin{center}
\begin{figure}[htbp]
\centerline{\includegraphics[width=10cm,height=8cm]{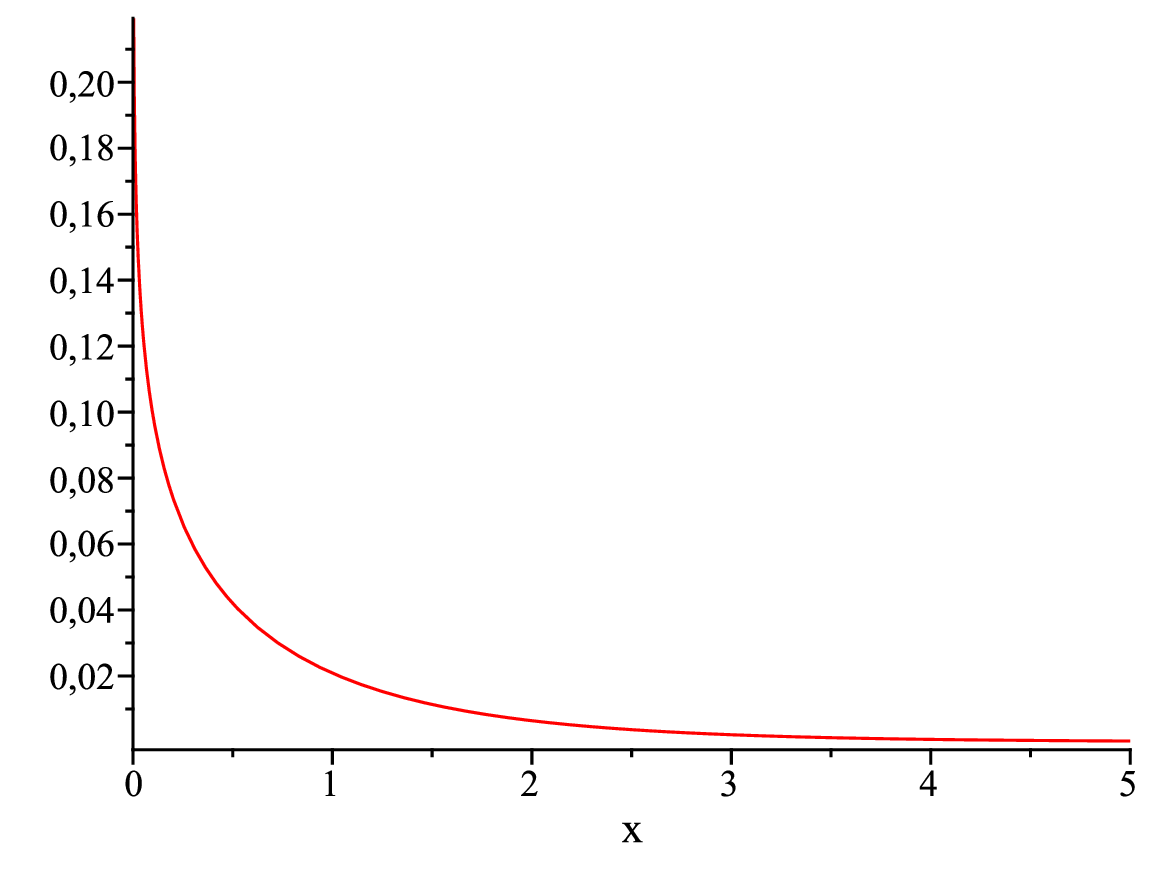}}
\caption{\it The shape of density $p_h(\bold x)$ on the interval $\Vert\bold x\Vert\in (0, 5]$
             (for $\lambda = 1, \mu=2, c=3$)} 
\label{FigPollutionDens}
\end{figure}
\end{center} 

Numerical calculations related to the values of density (\ref{eq2.0}) on the interval 
$\Vert\bold x\Vert\in (0,4]$ are also presented in Table 1 below. 

\begin{center}
\begin{tabular}{|r|r||r|r|}
\hline
 $\Vert\bold x\Vert$ & $p_l(\bold x)$ & $\Vert\bold x\Vert$ & $p_l(\bold x)$ \\
\hline\hline
  0.2 & 0.074088 & 2.2 & 0.005156 \\
\hline
  0.4 & 0.049587 & 2.4 & 0.004146 \\
\hline
  0.6 & 0.036023 & 2.6 & 0.003342 \\
\hline
  0.8 & 0.027129 & 2.8 & 0.002699 \\
\hline
  1.0 & 0.020854 & 3.0 & 0.002184 \\
\hline
  1.2 & 0.016243 & 3.2 & 0.001770 \\
\hline
  1.4 & 0.012771 & 3.4 & 0.001436 \\
\hline
  1.6 & 0.010110 & 3.6 & 0.001167 \\
\hline
  1.8 & 0.008046 & 3.8 & 0.000949 \\
\hline
  2.0 & 0.006430 & 4.0 & 0.000772 \\
\hline 
\end{tabular} 
\end{center}

\begin{center}
{\bf Table 1:} {\it Values of density $p_h(\bold x)$ on the interval $\Vert\bold x\Vert\in (0, 4]$ 
(for $\lambda=1, \mu=2, c=3$)} 
\end{center}

We see that, as should be expected, density (\ref{eq2.0}) takes its maximum value at the origin 
(that is, the level of pollution is maximal in the neighbourhood of the source), while it is decreasing  nonlinearly as the distance from the source increases. For example, on the unit circumference  
$\Vert\bold x\Vert = 1$, density (\ref{eq2.0}) takes the value: 
$p_h(\bold x)|_{\Vert\bold x\Vert = 1} \approx 0.020854$. We also see that for $\Vert\bold x\Vert > 3$ 
the density values become very close to zero, meaning that the vast majority of the pollutants will be concentrated in the circle $\Vert\bold x\Vert < 3$ around the source over time. 

The stationary density (\ref{eq2.0}) can also be used to approximately estimate the level of pollution 
in a given planar area for sufficiently large time values. Function (\ref{eq2.0}) shows how the density of pollution behaves, as the distance $\Vert\bold x\Vert$ from the source grows. To find the concentration of pollution in some planar area, say, in the circle $C_r = \{\bold x\in\Bbb R^2 : \; \Vert\bold x\Vert<r\}$ 
of some radius $r>0$, one needs to integrate density (\ref{eq2.0}) in $C_r$. The obtained value 
$k_r$ is less than 1 and yields the share of the total mass $M(t)$ of polluting substance emitting by a given time $t$ and settled in this circle. The product $k_r M(t)$ is the total concentration of polluting substance in $C_r$. Since the total mass $M(t)$ increases, as time $t$ grows, the concentration of 
polluting substance in $C_r$ increases too.

This interpretation allows us to predict the level of pollution in a given planar area on long time intervals. Moreover, knowledge of the distribution of polluting substance on the surface of a planar area enables us to pose a respective initial-value and/or boundary problem for modeling the process of percolating the polluting substance into lower soil strata, as described in \cite{stag}. In this case, the distribution of polluting substance yields the first initial and/or boundary condition.

\section{Light-particle diffusion model} 

Consider now a diffusion model when the polluting particles have relatively light random weight. 
In this case only a small their portion can settle near the source, while the overwhelming 
their majority falls at some distance from the source and the number of fallen particles is monotonously  and nonlinearly decreasing, as this distance increases. Such interpretation implies that the lifetime 
should be a positive random variable, whose density is an asymmetric bell-shaped function with a maximal   value at a point located at some distance from the source, and whose left slope is steeper than the 
right one.   

These are exactly the properties of the density of gamma-distribution with parameters 
$(\mu,\alpha), \; \mu>0, \; \alpha>2$, given by the formula: 
\begin{equation}\label{LifetimeDensity2}
q(t) = \left\{ 
\aligned  \frac{\mu^{\alpha}}{\Gamma(\alpha)} \; t^{\alpha-1} \; e^{-\mu t}, \quad & \text{if} \; t>0,\\ 
          0, \qquad \; & \text{if} \; t\le 0 ,  
\endaligned\right.
\qquad \mu>0 , \quad \alpha>2. 
\end{equation}
Note that for $\alpha=1$, density (\ref{LifetimeDensity2}) turns into exponential density 
(\ref{LifetimeDensity1}) related to the heavy-particle model. That is why we consider the case 
$\alpha>2$ only. The case $1<\alpha<2$, corresponding to the semi-heavy particles can also be examined 
in the same manner. Density (\ref{LifetimeDensity2}) takes maximal value at the point $t=(\alpha-1)/\mu$ 
and, therefore, by varying the value of parameter $\alpha$, we can shift this point changing the shape of 
the density and considering other models (semi-light, super-light weights, etc).  

According to general formula (\ref{PollutionDensityGeneral}), the density of pollution at the 
distance $\Vert\bold x\Vert$ from the source is given by the formula: 
\begin{equation}\label{PollutionDensity2}
p_l(\bold x)  = \frac{\lambda \; \mu^{\alpha}}{2\pi c \; \Gamma(\alpha)} \; 
\int_{\Vert\bold x\Vert/c}^{\infty} e^{-(\lambda+\mu)\tau} \; \tau^{\alpha-1} \; 
\frac{\exp \left(\frac{\lambda}{c} \sqrt{c^2\tau^2 - \Vert\bold x\Vert^2} \right)}{\sqrt{c^2\tau^2 
- \Vert\bold x\Vert^2}} \; d\tau , 
\end{equation}
where $\lambda$ is the intensity of changes of direction and $(\mu, \alpha)$ are the parameters of the 
gamma-distribution of the particle's lifetime with density (\ref{LifetimeDensity2}). 

\bigskip 

{\bf Theorem 2.} {\it The stationary probability density $p_l(\bold x)$ of the symmetric 
planar Markov random flight with $(\mu,\alpha)$-gamma distributed random lifetime 
(\ref{LifetimeDensity2}), is given by the formula:} 

\begin{equation}\label{eq3.0} 
\aligned 
p_l(\bold x) & = \frac{\lambda \; \mu^{\alpha}}{4\pi c^2 \; \Gamma(\alpha)} \; \sum_{k=0}^{\infty} 
\frac{\lambda^k}{k!}  \biggl[ \left( \frac{\Vert\bold x\Vert}{c} \right)^{\alpha+k-1} 
B\left( \frac{k+1}{2}, - \frac{\alpha+k-1}{2} \right) \\ 
& \hskip 4cm \times \; _1F_2\left( \frac{\alpha}{2};  \frac{1}{2}, \frac{\alpha+k+1}{2}; (\lambda+\mu)^2 
\; \frac{\Vert\bold x\Vert^2}{4c^2} \right) \\ 
& \hskip 3cm - (\lambda+\mu) \left( \frac{\Vert\bold x\Vert}{c} \right)^{\alpha+k} 
B\left( \frac{k+1}{2}, - \frac{\alpha+k}{2} \right) \\ 
& \hskip 4cm \times \; _1F_2\left( \frac{1+\alpha}{2};  \frac{3}{2}, \frac{\alpha+k}{2}+1; (\lambda+\mu)^2 
\; \frac{\Vert\bold x\Vert^2}{4c^2} \right) \\ 
& \hskip 3cm + 2^{2-\alpha-k} \; (\lambda+\mu)^{1-\alpha-k} \; \Gamma(\alpha+k-1) \\  
& \hskip 4cm \times \;  _1F_2\left( -\frac{k-1}{2}; -\frac{\alpha+k-3}{2}, -\frac{2\alpha+k+1}{4}; 
(\lambda+\mu)^2 \; \frac{\Vert\bold x\Vert^2}{4c^2} \right) \biggr] ,  
\endaligned
\end{equation}
{\it where} 
\begin{equation}\label{1F2}
_1F_2(\xi; \eta,\zeta; z) = \sum_{k=0}^{\infty} \frac{(\xi)_k}{(\eta)_k \; (\zeta)_k} \; \frac{z^k}{k!} 
\end{equation}
{\it is the hypergeometric function and $B(x,y)$ is the beta-function.}

\begin{proof} 

By expanding the exponential function in the fraction of the integrand in (\ref{PollutionDensity2}), 
we obtain: 
$$\aligned
& p_l(\bold x) \\ 
& = \frac{\lambda \; \mu^{\alpha}}{2\pi c \; \Gamma(\alpha)} \; 
\int_{\Vert\bold x\Vert/c}^{\infty} e^{-(\lambda+\mu)\tau} \; \tau^{\alpha-1} \; 
\frac{\exp \left(\frac{\lambda}{c} \sqrt{c^2\tau^2 - \Vert\bold x\Vert^2} \right)}{\sqrt{c^2\tau^2 
- \Vert\bold x\Vert^2}} \; d\tau \\ 
& = \frac{\lambda \; \mu^{\alpha}}{2\pi c \; \Gamma(\alpha)} \; \int_{\Vert\bold x\Vert/c}^{\infty} 
e^{-(\lambda+\mu) \tau} \; \tau^{\alpha-1} \; 
\sum_{k=0}^{\infty} \frac{1}{k!} \; \left( \frac{\lambda}{c} \right)^k \; \left( \sqrt{c^2\tau^2 - 
\Vert\bold x\Vert^2} \right)^{k-1} \; d\tau \\ 
& = \frac{\lambda \; \mu^{\alpha}}{2\pi c^2 \; \Gamma(\alpha)} \sum_{k=0}^{\infty} \frac{\lambda^k}{k!} \; 
\int_{\Vert\bold x\Vert/c}^{\infty} \tau^{\alpha-1} \; e^{-(\lambda+\mu) \tau} \; \left( \tau^2 - 
\frac{\Vert\bold x\Vert^2}{c^2} \right)^{(k-1)/2} \; d\tau \\ 
& \qquad (\text{see \cite[item 2.3.7, formula 3]{pbm2}}) \\ 
& = \frac{\lambda \; \mu^{\alpha}}{4\pi c^2 \; \Gamma(\alpha)} \; \sum_{k=0}^{\infty} \frac{\lambda^k}{k!}  
\biggl[ \left( \frac{\Vert\bold x\Vert}{c} \right)^{\alpha+k-1} 
B\left( \frac{k+1}{2}, - \frac{\alpha+k-1}{2} \right) \\ 
& \hskip 4cm \times \; _1F_2\left( \frac{\alpha}{2};  \frac{1}{2}, \frac{\alpha+k+1}{2}; (\lambda+\mu)^2 
\; \frac{\Vert\bold x\Vert^2}{4c^2} \right) \\ 
& \hskip 3cm - (\lambda+\mu) \left( \frac{\Vert\bold x\Vert}{c} \right)^{\alpha+k} 
B\left( \frac{k+1}{2}, - \frac{\alpha+k}{2} \right) \\ 
& \hskip 4cm \times \; _1F_2\left( \frac{1+\alpha}{2};  \frac{3}{2}, \frac{\alpha+k}{2}+1; (\lambda+\mu)^2 
\; \frac{\Vert\bold x\Vert^2}{4c^2} \right) \\ 
& \hskip 3cm + 2^{2-\alpha-k} \; (\lambda+\mu)^{1-\alpha-k} \; \Gamma(\alpha+k-1) \\  
& \hskip 4cm \times \;  _1F_2\left( -\frac{k-1}{2}; -\frac{\alpha+k-3}{2}, -\frac{2\alpha+k+1}{4}; 
(\lambda+\mu)^2 \; \frac{\Vert\bold x\Vert^2}{4c^2} \right) \biggr] , 
\endaligned$$
where $_1F_2(\xi; \eta,\zeta; z)$ is the hypergeometric function defined by (\ref{1F2}) and $B(x,y)$ 
is the beta-function, also called the Euler integral of first kind. 

The theorem is proved.  
\end{proof}

Formula (\ref{eq3.0}) seems very complicated for practical numerical calculations because it has the form 
of a functional series composed of hypergeometric functions with variable coefficients and this fact implies fairly big computational difficulties. That is why it is more convenient to make numerical calculations  directly by means of formula (\ref{PollutionDensity2}). The results of such numerical calculations on the 
interval $\Vert\bold x\Vert\in [0, 4)$ (for the values of parameters $\lambda=1, \mu=2, c=2, \alpha=5$) 
are presented in Table 2 below. 

\begin{center}
\begin{tabular}{|r|r||r|r|}
\hline
 $\Vert\bold x\Vert$ & $p_l(\bold x)$ & $\Vert\bold x\Vert$ & $p_l(\bold x)$ \\
\hline\hline
  0.0 & 0.019894 & 2.0 & 0.016472 \\
\hline
  0.2 & 0.019896 & 2.2 & 0.015442 \\
\hline
  0.4 & 0.019901 & 2.4 & 0.014339 \\
\hline
  0.6 & 0.019887 & 2.6 & 0.013196 \\
\hline
  0.8 & 0.019814 & 2.8 & 0.012040 \\
\hline
  1.0 & 0.019636 & 3.0 & 0.010897 \\
\hline
  1.2 & 0.019317 & 3.2 & 0.009789 \\
\hline
  1.4 & 0.018837 & 3.4 & 0.008731 \\
\hline
  1.6 & 0.018194 & 3.6 & 0.007736 \\
\hline
  1.8 & 0.017399 & 3.8 & 0.006811 \\
\hline 
\end{tabular} 
\end{center}

\begin{center}
{\bf Table 2:} {\it Values of density $p_l(\bold x)$ on the interval $\Vert\bold x\Vert\in [0, 4)$ \\
(for $\lambda=1, \mu=2, c=2, \alpha=5$)} 
\end{center}

We see that, on the interval $\Vert\bold x\Vert\in [0.0, \; 0.4)$, density $p_l(\bold x)$ increases and 
takes its maximal value at the point $\Vert\bold x\Vert = 0.4$. Then, on the interval 
$\Vert\bold x\Vert\in (0.4, \; 4.0)$, this density becomes a monotonously decreasing function. As  
should be expected, such behaviour demonstrates the fact that the overwhelming majority of light particles  settle at some distance from the source. As the value of parameter $\alpha$ increases, the point of maximum density shifts in the positive direction, which corresponds to the model of more light particles.

\section{Some remarks on the asymmetric case} 

In the previous sections, we have presented a diffusion model of surface soil pollution based on the symmetric planar finite-velocity stochastic motion of the particles whose lifetime, that is, the time of active evolution, is a positive random variable with given distribution. The assumption of symmetry of such a  stochastic motion implies that each particle can take its directions (the initial and each new one) 
according to the uniform probability law. In practice, this means that the wind rose at the source point 
has the uniform nature. 

In reality, however, the wind rose can have a dominant direction and the diffusion area in this case is not a circle, but an asymmetric oval-shaped planar region, elongated from the source point in this dominant direction. Nevertheless, the same approach can also be used in this asymmetric case. For example, the asymmetric wind rose can be modeled by the von Mises distribution on the unit circumference 
$S_1$ with the two-dimensional density: 
\begin{equation}\label{ADFgauss1}
\chi_k(\bold x) = \frac{1}{2\pi \; I_0(k)} \; \exp\left( \frac{k x_1}{\Vert\bold x\Vert} \right) \; 
\delta(1-\Vert\bold x\Vert^2) , 
\end{equation}
$$\bold x =(x_1,x_2)\in\Bbb R^2, \qquad \Vert\bold x\Vert=\sqrt{x_1^2+x_2^2} \qquad k\in\Bbb R^1,$$
where $I_0(z)$ is the modified Bessel function of order zero. Formula (\ref{ADFgauss1}) determines the 
one-parametric family of densities $\left\{ \chi_k(\bold x), \; k\in\Bbb R^1 \right\}$ and, for any fixed real $k\in\Bbb R^1$, the density $\chi_k(\bold x)$ is absolutely continuous and uniformly bounded on 
$S_1$. If $k=0$, then formula (\ref{ADFgauss1}) yields the density of the uniform distribution 
on $S_1$, while for $k\neq 0$ it produces non-uniform densities.  

In the unit polar coordinates $x_1=\cos\theta, \; x_2=\sin\theta$, two-dimensional density 
(\ref{ADFgauss1}) takes the form of the circular Gaussian law:  
$$\chi_k(\theta) = \frac{\exp(k \cos\theta)}{2\pi \; I_0(k)} , \qquad \theta\in [-\pi, \pi), 
\quad k\in\Bbb R^1.$$   
In each particular case, the value of parameter $\theta$ should be taken in such a way that to get the respective  dominated direction. 

The asymmetric planar stochastic motion whose random directions are taken according to 
(\ref{ADFgauss1}), has the transition probability density given by (see \cite[formula (4.12.34)]{kol1}):
\begin{equation}\label{ADFgauss6}
p(\bold x,t) = \sum_{n=0}^{\infty} \lambda^n \left( \frac{1}{2\pi c \; I_0(k)} \right)^{n+1} 
\left[ \frac{e^{-\lambda t}}{t} \; \exp\left( \frac{k x_1}{\Vert\bold x\Vert} \right) \; 
\delta(c^2 t^2-\Vert\bold x\Vert^2) \right]^{\overset{\bold x}{\ast} \overset{t}{\ast}(n+1)} , 
\end{equation} 
$$\bold x =(x_1,x_2)\in\Bbb R^2, \qquad \Vert\bold x\Vert=\sqrt{x_1^2+x_2^2} \qquad k\in\Bbb R^1,$$ 
where the symbol ${\overset{\bold x}{\ast} \overset{t}{\ast}(n+1)}$ means the $(n+1)$-multiple double convolution with respect to variables $\bold x$ and $t$. 

Then, similarly as above, the stationary density of pollution in the asymmetric case is given by 
the general formula: 
\begin{equation}\label{PollutionDensityAsym}
p_{asym}(\bold x) =  \sum_{n=0}^{\infty} \lambda^n \left( \frac{1}{2\pi c \; I_0(k)} \right)^{n+1} 
\int_{\Vert\bold x\Vert/c}^{\infty} \left[ \frac{e^{-\lambda\tau}}{\tau} \; \exp\left( \frac{k x_1}{\Vert\bold x\Vert} \right) \; \delta(c^2\tau^2-\Vert\bold x\Vert^2) \right]^{\overset{\bold x}{\ast} 
\overset{\tau}{\ast}(n+1)} q(\tau) \; d\tau , 
\end{equation}
where, as above, $q(t)$ is the density of the particle's lifetime. 

Evaluation of density (\ref{PollutionDensityAsym}) is an extremely difficult, and even impracticable, 
analytical problem determined by the possibility of calculating the double convolutions of the singular  function in square brackets of the integrand. There is, however, a hope to numerically calculate at least 
a few first terms of this series.

\section{Conclusions} 

In this article, a diffusion model of surface soil pollution based on planar symmetric Markov 
random flight, is presented. The basic feature of the model is the assumption that the pollution process 
is carried out by the particles moving chaotically with finite speed in the plane and whose time of active evolution, called random lifetime, is a random variable with given distribution. The key point in applying this model is the choice of an appropriate distribution of the particles' lifetime, which crucially depends on their weights. 

In practice, to choose a suitable lifetime distribution, it is necessary to conduct the statistical analysis of the particles' weights, construct a histogram and, based on its shape, select a suitable empirical density with the appropriate parameters (exponential, gamma or some other). Such a histogram and density will show what type of particles are involved in the pollution process (heavy, light, semi-light, super-heavy, etc.). The parameters $c$ and $\lambda$ can be selected in such a way as to ensure the maximum closeness of the empirical density to the histogram. Substituting then the empirical density into integral 
(\ref{PollutionDensityGeneral}) in the symmetric case, or into series (\ref{PollutionDensityAsym}) in the asymmetric case, we obtain a general expression for the stationary pollution density.  

Evaluation of such expression is a very difficult analytical problem that, apparently, can be calculated numerically only.

\bigskip 

{\bf Acknowledgements.} This research was supported by the National Agency for Research and Development 
(ANCD) of Moldova in the framework of the project No. 24.80012.5007.25SE.  

\bigskip 

{\bf Declaration.} The author declares no potential conflicts of interest with respect to the research, authorship, and/or publication of this article. The author has no data availability to share.

\end{document}